\documentclass[12pt]{article}
\usepackage{amsfonts}
\usepackage{amsmath,amssymb,amsfonts}
\usepackage{color}
\usepackage{hyperref}

\def\squarebox#1{\hbox to #1{\hfill\vbox to #1{\vfill}}}
\newcommand{\qed}{\hspace*{\fill}
\vbox{\hrule\hbox{\vrule\squarebox{.667em}\vrule}\hrule}\smallskip}
\newtheorem{teorema}{Theorem}[section]
\newtheorem{lema}[teorema]{Lemma}
\newtheorem{corolario}[teorema]{Corollary}
\newtheorem{proposicao}[teorema]{Proposition}

\newenvironment{profe}{\noindent {\bf Proof:}}{\hfill $\qed $ \newline}
\begin{document}

\title{Controllability of control systems simple Lie groups and the topology
of flag manifolds}
\author{Ariane Luzia dos Santos\thanks{Supported by CNPq grant 
n$^{\mathrm{o} }$ 142536/2008-3}\thanks{%
Address: Imecc - Unicamp, Departamento de Matem\'{a}tica. Rua S\'{e}rgio
Buarque de Holanda, 651, Cidade Universit\'{a}ria Zeferino Vaz. 13083-859
Campinas S\~{a}o Paulo, Brasil. e-mail: ra069520@ime.unicamp.br arianeluzsantos@gmail.com}
 \and Luiz A. B. San Martin\thanks{%
Supported by CNPq grant n$^{\mathrm{o}}$ 303755/2009-1 and FAPESP grant n$^{%
\mathrm{o}}$ 07/06896-5}\thanks{%
Address: Imecc - Unicamp, Departamento de Matem\'{a}tica. Rua S\'{e}rgio
Buarque de Holanda, 651, Cidade Universit\'{a}ria Zeferino Vaz. 13083-859
Campinas S\~{a}o Paulo, Brasil. e-mail: smartin@ime.unicamp.br}}
\date{}
\maketitle

\begin{abstract}
Let $S$ be subsemigroup with nonempty interior of a complex simple Lie group 
$G$. It is proved that $S=G$ if $S$ contains a subgroup $G\left( \alpha
\right) \approx \mathrm{Sl}\left( 2,\mathbb{C}\right) $ 
generated by the $\exp \mathfrak{g}_{\pm \alpha }$, where $\mathfrak{g}%
_{\alpha }$ is the root space of the root $\alpha $. The proof uses the
fact, proved before, that the invariant control set of $S$ is contractible
in some flag manifold if $S$ is proper, and exploits the fact that several
orbits of $G\left( \alpha \right) $ are $2$-spheres not null homotopic. The
result is applied to revisit a controllability theorem 
and get some improvements. 
\end{abstract}

%TCIMACRO{\TeXButton{noindent}{\noindent}}%
%BeginExpansion
\noindent%
%EndExpansion

%TCIMACRO{\TeXButton{noindent}{\noindent}}%
%BeginExpansion
\noindent%
%EndExpansion

\section{Introduction}

In  this paper we use a method to study controllability of bilinear control
systems and invariant control systems on (semi-)simple Lie groups that
relies on the (algebraic) topology of flag manifolds.

The method is based on the geometry of invariant control sets on flag
manifolds, as described initially in \cite{sminv}, \cite{SMT} and \cite%
{smsurv}, and further developed in \cite{SM},  \cite{smbrsimplt}, \cite%
{smmax}, \cite{smord}, \cite{smroconex}, \cite{smsanhom}.

The part of this  geometry to be applied here states that if $S\subset G$ is
a semigroup with nonempty interior  then there exists some flag manifold of $%
G$, say $\mathbb{F}_{\Theta }$, such that the unique invariant control set $%
C_{\Theta }\subset \mathbb{F}_{\Theta }$, for the action of $S$  on $\mathbb{%
F}_{\Theta }$, is contained in a subset $\mathcal{E}\subset \mathbb{F}%
_{\Theta }$, which is homeomorphic to an Euclidian space $\mathbb{R}^{N}$
(cf. Theorem \ref{teoregre} below). ($\mathcal{E}\approx \mathbb{R}^{N}$ is
an open Bruhat cell of $\mathbb{F}_{\Theta }$.)

This implies, for instance, that any closed curve $\gamma $ contained in $%
C_{\Theta }$ is homotopic (in $\mathbb{F}_{\Theta }$) to a point, and hence
represents a trivial element of the fundamental group $\pi _{1}\left( 
\mathbb{F}_{\Theta }\right) $. Analogously, any higher dimensional sphere $%
S^{n}\subset C_{\Theta }$ represents the identity of the homotopy group $\pi
_{n}\left( \mathbb{F}_{\Theta }\right) $.

Therefore one can achieve to prove $S=G$  by showing that the invariant
control sets in all flag manifolds are topologically non trivial e.g.
contains a curve or a sphere not homotopic (in the flag manifold) to a
point. In particular one gets controllability (in $G$) of an invariant
control system on $G$ by applying this method to the semigroup of control $S$%
, as soon as the control system satisfies the Lie algebra rank condition.

In this paper we give sufficient conditions for controllability (under the
Lie algebra rank condition) taking advantage of the fact that some subgroups
of $G$ have orbits on the flag manifolds that are homeomorphic to spheres
but not homotopic to a point. Specifically, we consider subgroups $G\left(
\alpha \right) \subset G$ where $\alpha $ is a root of the Lie algebra $%
\mathfrak{g}$ of $G$ and $G\left( \alpha \right) $ is generated by (the
exponentials of)  the root spaces $\mathfrak{g}_{\alpha }$ and $\mathfrak{g}%
_{-\alpha }$.

Then our main result (see Theorem \ref{teocomplexo}  below) says that in a
complex Lie group $G$ the semigroup $S=G$ if $\mathrm{int}S\neq \emptyset $
and $G\left( \alpha \right) \subset S$ The technique of proof consist in i)
checking that several orbits $G\left( \alpha \right) \cdot x$ are $2$%
-spheres not homotopic to a point; and ii) some of these orbits are
contained in the unique invariant control set $C_{\Theta }$ of $S$. If this
is done in any flag manifold $\mathbb{F}_{\Theta }$ then no $C_{\Theta }$ is
contained in a contractible subset, and $S$ must be $G$.

Our source of inspiration to think in the group $G\left( \alpha \right) $ is
a series of papers started with Jurdjevic-Kupka \cite{JK}, \cite{JK1},
followed by several papers (see Gauthier-Kupka-Sallet \cite{GKS} and
references therein), and culminating with the final result of El
Assoudi-Gauthier-Kupka \cite{EGK}. One of the main issues in these papers is
that the semigroup of control $S$ contains a regular element as well as $%
G\left( \mu \right) $ when $\mu $ is the highest root.

Thus our Theorem \ref{teocomplexo} provides an alternate proof of the main
theorem of \cite{EGK}, when the group $G$ is simple and complex. Actually,
we improve that result for these groups. This is because our result is for
an arbitrary root $\alpha $, and not just the highest one. We state this
improvement in Theorem \ref{teocontrolreal}.

We work with simple groups to avoid to take all the time the decomposition
into simple components, which can be done in the standard fashion, and is
left to the reader.

Similar results can be obtained for  real simple groups although the
topology of their flag manifolds is trickier. We leave to a forthcoming
paper the case of  the so-called normal real forms  where all the roots have
multiplicity one, and hence  the orbits $G\left( \alpha \right) \cdot
b_{\Theta }$ have dimension zero or one. In this case  we must look at the
fundamental groups $\pi _{1}\left( \mathbb{F}_{\Theta }\right) $.

\section{Semigroups and flag manifolds}

For  a complex semi-simple Lie group $G$ with Lie algebra $\mathfrak{g}$ we
use the following notation:

\begin{itemize}
\item $\mathfrak{h}$ is a Cartan subalgebra, whose set of roots is denoted
by $\Pi $. $\Pi ^{+}$ is a set of positive roots with 
\begin{equation*}
\Sigma =\{\alpha _{1},\ldots ,\alpha _{l}\}\subset \Pi ^{+} 
\end{equation*}
standing for the corresponding simple system of roots. We have $\Pi =\Pi ^{+}%
\dot{\cup}\left( -\Pi ^{+}\right) $ and any $\alpha \in \Pi ^{+}$, is a
linear combination $\alpha =n_{1}\alpha _{1}+\cdots +n_{l}\alpha _{l}$ with $%
n_{i}\geq 0$ integers. The support of $\alpha $, $\mathrm{supp}\alpha $ is
the subset of $\Sigma $ where $n_{i}>0$.

\item The Cartan-Killing form of $\mathfrak{g}$ is denoted by $\langle \cdot
,\cdot \rangle $. If $\alpha \in \mathfrak{h}^{*}$ then $H_{\alpha }\in 
\mathfrak{h}$ is defined by $\alpha \left( \cdot \right) =\langle H_{\alpha
},\cdot \rangle $, and $\langle \alpha ,\beta \rangle =\langle H_{\alpha
},H_{\beta }\rangle $. The subspace spanned over $\mathbb{R}$ by $H_{\alpha }
$, $\alpha \in \Pi $, is denoted by $\mathfrak{h}_{\mathbb{R}}$. We have $%
\mathfrak{h}=\mathfrak{h}_{\mathbb{R}}+i\mathfrak{h}_{\mathbb{R}}$.

\item We write 
\begin{equation*}
\mathfrak{h}_{\mathbb{R}}^{+}=\{H\in \mathfrak{h}_{\mathbb{R}}:\forall
\alpha \in \Pi ^{+},\,\alpha \left( H\right) >0\} 
\end{equation*}
for the Weyl chamber defined by $\Pi ^{+}$.

\item The root space of a root $\alpha $ is 
\begin{equation*}
\mathfrak{g}_{\alpha }=\{X\in \mathfrak{g}:\forall H\in \mathfrak{h}%
,\,[H,X]=\alpha \left( H\right) X\}. 
\end{equation*}
It is known that $\dim _{\mathbb{C}}\mathfrak{g}_{\alpha }=1$.

\item For a root $\alpha $, $\mathfrak{g}\left( \alpha \right) $ is the
subalgebra generated by $\mathfrak{g}_{\alpha }$ and $\mathfrak{g}_{-\alpha }
$. Then 
\begin{equation*}
\mathfrak{g}\left( \alpha \right) =\mathrm{span}_{\mathbb{C}}\{H\}\oplus 
\mathfrak{g}_{\alpha }\oplus \mathfrak{g}_{-\alpha }\approx \mathfrak{sl}%
\left( 2,\mathbb{C}\right) .
\end{equation*}%
$G\left( \alpha \right) $ is the connected Lie subgroup with Lie algebra $%
\mathfrak{g}\left( \alpha \right) $, which is isomorphic to $\mathrm{Sl}%
\left( 2,\mathbb{C}\right) /D$, where $D$ is a discrete central subgroup
(because $\mathrm{Sl}\left( 2,\mathbb{C}\right) $ is simply connected).

\item $\mathfrak{u}$ is a compact real form of $\mathfrak{g}$ and $U=\langle
\exp \mathfrak{u}\rangle $ is the connected Lie subgroup with Lie algebra $%
\mathfrak{u}$. It is known that $U$ is compact semi-simple, and maximal
compact in $G$.

\item $\mathcal{W}$ is the Weyl group. Either $\mathcal{W}$ is the group
generated by the reflections $r_{\alpha }$, $\alpha \in \Pi $, $r_{\alpha
}\left( \beta \right) =\beta -\frac{2\langle \alpha ,\beta \rangle }{\langle
\alpha ,\alpha \rangle }\alpha $, or $\mathcal{W}=\mathrm{Norm}_{U}\left( 
\mathfrak{h}\right) /T$ where $T$ is the torus $U\cap \exp \mathfrak{h}$ and 
$\mathrm{Norm}_{U}\left( \mathfrak{h}\right) =\{g\in U:\mathrm{Ad}\left(
g\right) \mathfrak{h}\subset \mathfrak{h}\}$ is the normalizer of $\mathfrak{%
h}$ in $U$.

\item $\mathfrak{n}^{+}=\sum_{\alpha \in \Pi ^{+}}\mathfrak{g}_{\alpha }$
and $\mathfrak{n}^{-}=\sum_{\alpha \in \Pi ^{+}}\mathfrak{g}_{-\alpha }$

\item Given the data $\mathfrak{h}$ and $\Pi ^{+}$ (or $\Sigma $) there is
the Borel subalgebra (minimal parabolic) $\mathfrak{p}=\mathfrak{h}\oplus 
\mathfrak{n}^{+}$. A subset $\Theta \subset \Sigma $ defines the standard
parabolic subalgebra by 
\begin{equation*}
\mathfrak{p}_{\Theta }=\mathfrak{p}+\sum_{\alpha \in \langle \Theta \rangle }%
\mathfrak{g}_{\alpha } 
\end{equation*}
where $\langle \Theta \rangle =\{\alpha \in \Pi :$ $\mathrm{supp}\alpha
\subset \Theta $ or $\mathrm{supp}\left( -\alpha \right) \subset \Theta \}$
is the set of roots spanned by $\Theta $. ($\mathfrak{p}_{\emptyset }=%
\mathfrak{p}$.)

\item For $\Theta \subset \Sigma $, $P_{\Theta }$ is the parabolic subgroup
with Lie algebra $\mathfrak{p}_{\Theta }$, which is the normalizer of $%
\mathfrak{p}_{\Theta }$: 
\begin{equation*}
P_{\Theta }=\mathrm{Norm}_{U}\left( \mathfrak{p}_{\Theta }\right) =\{g\in U:%
\mathrm{Ad}\left( g\right) \mathfrak{p}_{\Theta }\subset \mathfrak{p}%
_{\Theta }\}.
\end{equation*}

\item The flag manifold $\mathbb{F}_{\Theta }=G/P_{\Theta }$, which is
independent of the specific group $G$ with Lie algebra $\mathfrak{g}$. The
origin of $G/P_{\Theta }$ (the coset $1\cdot P_{\Theta }$) is denoted by $%
b_{\Theta }$.
\end{itemize}

Now let $S\subset G$ be a subsemigroup with $\mathrm{int}S\neq \emptyset $.
We recall here some results of \cite{sminv}, \cite{smsurv} and \cite{SMT}
that are on the basis of our topological approach to controllability in $G$.

We let $S$ act on a flag manifold $\mathbb{F}_{\Theta }$, by restriction of
the action of $G$. An invariant control set for $S$ in $\mathbb{F}_{\Theta }$
is a subset $C\subset \mathbb{F}_{\Theta }$ such that $\mathrm{cl}\left(
Sx\right) =C $ for every $x\in C$, where $Sx=\{gx\in \mathbb{F}_{\Theta
}:g\in S\}$. Since $\mathrm{int}S\neq \emptyset $ such a set is closed, has
nonempty interior and is in fact invariant ($gx\in C$ if $g\in S$ and $x\in C
$).

\begin{lema}
\label{lemuniqueics}(See \cite{sminv}.) In any flag manifold $\mathbb{F}%
_{\Theta }$ there is a unique invariant control set for $S$, denoted by $%
C_{\Theta }$.
\end{lema}

To state the geometrical property of $C_{\Theta }$ to be used later we
discuss the dynamics of the vector fields $\widetilde{H}$ on a flag manifold 
$\mathbb{F}_{\Theta }$ whose flow is $\exp tH$, with $H$ in the closure $%
\mathrm{cl}\mathfrak{h}_{\mathbb{R}}^{+}$ of Weyl chamber $\mathfrak{h}_{%
\mathbb{R}}^{+}$. It is known that $\widetilde{H}$ is a gradient vector
field with respect to some Riemmannian metric on $\mathbb{F}_{\Theta }$ (see
Duistermat-Kolk-Varadarajan \cite{dkv} and Ferraiol-Patr\~{a}o-Seco \cite%
{fepase}).

Hence the orbits of $\widetilde{H}$ are either fixed points or trajectories
flowing between fixed point sets. Moreover, $\widetilde{H}$ has a unique
attractor fixed point set, say $\mathrm{att}_{\Theta }\left( H\right) $,
that has an open and dense stable manifold $\sigma _{\Theta }\left( H\right) 
$ (see \cite{dkv} and \cite{fepase}). This means that if $x\in \sigma
_{\Theta }\left( H\right) $ then its $\omega $-limit set $\omega \left(
x\right) $ is contained in $\mathrm{att}_{\Theta }\left( H\right) $. This
attractor has the following algebraic expressions 
\begin{equation*}
\mathrm{att}_{\Theta }\left( H\right) =Z_{H}\cdot b_{\Theta }=U_{H}\cdot
b_{\Theta } 
\end{equation*}
(see \cite{dkv} and \cite{fepase}). Here $Z_{H}=\{g\in G:\mathrm{Ad}\left(
g\right) H=H\}$ is the centralizer of $H$ in $G$ and $U_{H}=Z_{H}\cap U$ is
the centralizer in $U$. Its stable set $\sigma _{\Theta }\left( H\right) $
is also described algebraically by 
\begin{equation*}
\sigma _{\Theta }\left( H\right) =N_{H}^{-}Z_{H}\cdot b_{\Theta } 
\end{equation*}
where $N_{H}^{-}=\exp \mathfrak{n}_{H}^{-}$ and 
\begin{equation*}
\mathfrak{n}_{H}^{-}=\sum_{\gamma \left( H\right) <0}\mathfrak{g}_{\gamma }. 
\end{equation*}

In particular if $H$ is regular, that is, $H\in \mathfrak{h}_{\mathbb{R}}$
and $\alpha \left( H\right) >0$ for $\alpha \in \Pi ^{+}$ then $Z_{H}$
reduces to the Cartan subgroup $\exp \mathfrak{h}$, which fixes $b_{\Theta }$%
. Hence 
\begin{equation*}
\mathrm{att}_{\Theta }\left( H\right) =Z_{H}\cdot b_{\Theta }=\{b_{\Theta
}\}\qquad H\in \mathfrak{h}_{\mathbb{R}}.
\end{equation*}%
Actually, in the regular case the fixed points are isolated because $%
\widetilde{H}$ is the gradient of a Morse function, see \cite{dkv} and \cite%
{fepase}. Also, $\mathfrak{n}_{H}^{-}=\mathfrak{n}^{-}$ (notation as above)
and the stable set is $N^{-}\cdot b_{\Theta }$ (open Bruhat cell).

The following statement is a well known result from the Bruhat decomposition
of the flag manifolds (see \cite{dkv}, \cite{K}, \cite{W}).

\begin{proposicao}
In any flag manifold $\mathbb{F}_{\Theta }$ the open Bruhat cell $N^{-}\cdot
b_{\Theta }$ is diffeomorphic to an Euclidian space $\mathbb{R}^{d}$. (The
diffeomorphism is $X\in \mathfrak{n}_{\Theta }^{-}\mapsto \exp X\cdot
b_{\Theta } $, where $\mathfrak{n}_{\Theta }^{-}=\sum \{\mathfrak{g}_{\alpha
}:\alpha <0$ and $\alpha \notin \langle \Theta \rangle \}$.
\end{proposicao}

Put $h=\exp H$, $H\in \mathfrak{h}_{\mathbb{R}}^{+}$. It follows from the
gradient property of $\widetilde{H}$ that $\lim_{n\rightarrow +\infty
}h^{n}x=b_{\Theta }$ for any $x\in N^{-}\cdot b_{\Theta }$.

Now, we say that $g\in G$ is \textbf{regular real} if it is  a conjugate $%
g=aha^{-1}$ of $h=\exp H$, $H\in \mathfrak{h}_{\mathbb{R}}^{+}$ with $a\in G$%
. Then we write $\sigma _{\Theta }\left( g\right) =g\cdot \sigma _{\Theta
}\left( H\right) $ and call this the stable set of $g$ in $\mathbb{F}%
_{\Theta }$. (The reason for this name is clear: $g^{n}=\left(
aha^{-1}\right) ^{n}=ah^{n}a^{-1}$ and hence $g^{n}x\rightarrow gb_{\Theta }$
if $x\in \sigma _{\Theta }\left( g\right) $.)

The following lemma was used in \cite{sminv} to prove the above Lemma \ref%
{lemuniqueics}.

\begin{lema}
(See \cite{sminv}.) There exists regular real $g\in \mathrm{int}S$.
\end{lema}

Now we can state the next theorem from \cite{SMT}, which is in the basis of
our approach to controllability.

\begin{teorema}
\label{teoregre}Suppose that $S\neq G$. Then there exists a flag manifold $%
\mathbb{F}_{\Theta }$ such that the invariant control set $C_{\Theta
}\subset \sigma _{\Theta }\left( g\right) $ for every regular real $g\in 
\mathrm{int}S $.
\end{teorema}

\begin{corolario}
\label{coreuclid}If $S\neq G$ then $C_{\Theta }$ is contained in a subset $%
\mathcal{E}_{\Theta }\subset \mathbb{F}_{\Theta }$, which is diffeomorphic
to an Euclidian space.
\end{corolario}

%TCIMACRO{\TeXButton{vspace12pt}{\vspace{12pt}}}%
%BeginExpansion
\vspace{12pt}%
%EndExpansion

\noindent \textbf{Remark:} It can be proved that there exists a minimal $%
\Theta _{S}$ satisfying the condition of Theorem \ref{teoregre}. This $%
\Theta _{S}$ (or rather the flag manifold $\mathbb{F}_{\Theta _{S}}$) is
called the flag type of $S$ or the parabolic type of $S$ (because of the
parabolic subgroup $P_{\Theta _{S}}$). Several properties of $S$ are derived
from this flag type (e.g. the homotopy type of $S$ as in \cite{smsanhom} or
the connected components of $S$ as in \cite{smroconex}).

\section{Root spaces and semigroups}

Recall the notation $\mathfrak{g}\left( \alpha \right) =\mathrm{span}_{%
\mathbb{C}}\{H_{\alpha }\}\oplus \mathfrak{g}_{\alpha }\oplus \mathfrak{g}%
_{-\alpha }\approx \mathfrak{sl}\left( 2,\mathbb{C}\right) $ and $G\left(
\alpha \right) =\langle \exp \mathfrak{g}\left( \alpha \right) \rangle $.

\begin{teorema}
\label{teocomplexo}Let $G$ be a simple complex Lie group and $S\subset G$ a
semigroup with $\mathrm{int}S\neq \emptyset $. Then $S=G$ if there is a root 
$\alpha $ with $G\left( \alpha \right) \subset S$.
\end{teorema}

For the proof of this theorem we exploit the fact that $\mathfrak{g}\left(
\alpha \right) $ is isomorphic to $\mathfrak{sl}\left( 2,\mathbb{C}\right) $%
, and hence the unique flag manifold of $G\left( \alpha \right) $ is a $2$%
-sphere. With this in mind we prove that in any flag manifold $\mathbb{F}%
_{\Theta }$ of $G$ there are several $G\left( \alpha \right) $-orbits that
are $2$-spheres. Then we ensure that some of these orbits are contained in
the invariant control set $C_{\Theta }$ of $S$ in $\mathbb{F}_{\Theta }$.
Finally we use De Rham cohomology of $\mathbb{F}_{\Theta }$ to prove that
such orbits ($2$-spheres) are not homotopic to a point. Hence $C_{\Theta }$
is not contained in a contractible set and the theorem follows by Corollary %
\ref{coreuclid}.

The first step is the following lemma which reduces the proof to some
specific roots.

\begin{lema}
\label{lemactweyl}Let $\gamma $ and $\beta $ be roots such that $\beta
=w\gamma $ with $w\in \mathcal{W}$. If Theorem \ref{teocomplexo} holds for $%
\gamma $ then it is true for $\beta $ as well.
\end{lema}

\begin{profe}
Take a representative $\overline{w}$ of $w$ in the normalizer $M^{\ast }$ of 
$\mathfrak{h}$ in $U$. Then $\beta =w\gamma $ entails that $G\left( \beta
\right) =\overline{w}G\left( \gamma \right) \overline{w}^{-1}$. Now if $%
G\left( \beta \right) \subset S$ then $G\left( \gamma \right) $ is contained
in the $\overline{w}S\overline{w}^{-1}$ that has nonempty interior. Hence $%
\overline{w}S\overline{w}^{-1}=G$ implying that $S=G$.
\end{profe}

The assumption that $G$ is simple ensures the following fact about the
action Weyl group $\mathcal{W}$ on the set of roots $\Pi $: It is transitive
for the Dynkin diagrams $A_{l}$, $D_{l}$, $E_{6}$, $E_{7}$ and $E_{8}$, that
have only simple edges. For the other diagrams $B_{l}$, $C_{l}$, $F_{4}$ and 
$G_{2}$, there are two orbits which are given by the long and the short
roots, respectively.

On the other hand $\mathcal{W}$ acts transitively on the set of chambers.
Hence for any $\alpha \in \Pi $ there exists a unique root $\mu $ such that $%
H_{\mu }\in \mathrm{cl}\mathfrak{h}_{\mathbb{R}}^{+}$. Hence there is either
one or two roots with $H_{\mu }\in \mathrm{cl}\mathfrak{h}_{\mathbb{R}}^{+}$%
. By the above lemma it is enough to prove \ the theorem for these roots.

For the simply laced diagrams or for the long roots in the other diagrams
the root $\mu $ with $H_{\mu }\in \mathrm{cl}\mathfrak{h}_{\mathbb{R}}^{+}$
is the highest root. In any case we have the following property.

\begin{proposicao}
\label{propsuport}If $\mu $ is a root with $H_{\mu }\in \mathrm{cl}\mathfrak{%
h}_{\mathbb{R}}^{+}$ then $\mathrm{supp}\mu =\Sigma $.
\end{proposicao}

\begin{profe}
If $\{H_{1},\ldots ,H_{l}\}$ is the dual basis of $\Sigma =\{\alpha
_{1},\ldots ,\alpha _{l}\}$ then $\alpha =\alpha \left( H_{1}\right) \alpha
_{1}+\cdots +\alpha \left( H_{l}\right) \alpha _{l}$ for any $\alpha \in 
\mathfrak{h}^{\ast }$. The closure $\mathrm{cl}\mathfrak{h}_{\mathbb{R}}^{+}$
of the Weyl chamber is a polyhedral cone spanned by $H_{1},\ldots ,H_{l}$.
It is known that $\langle A,B\rangle >0$ for nonzero $A,B\in \mathrm{cl}%
\mathfrak{h}_{\mathbb{R}}^{+}$. Therefore, if $H_{\mu }\in \mathrm{cl}%
\mathfrak{h}_{\mathbb{R}}^{+}$ then the coefficientes $\mu \left(
H_{i}\right) $ satisfy $\mu \left( H_{i}\right) =$ $\langle H_{\mu
},H_{i}\rangle >0$, $i=1,\ldots ,l$, which means that $\mathrm{supp}\mu
=\Sigma $.
\end{profe}

From now on $\mu $ stands for one of the roots with $H_{\mu }\in \mathrm{cl}%
\mathfrak{h}_{\mathbb{R}}^{+}$.

For the next step we recall the notation of the last section where $\mathrm{%
att}_{\Theta }\left( H_{\mu }\right) $ is the attractor fixed point set of $%
\widetilde{H}_{\mu }\in \mathrm{cl}\mathfrak{h}_{\mathbb{R}}^{+}$ with $%
\sigma _{\Theta }\left( H_{\mu }\right) $ the stable set.

Now let $C_{\Theta }$ be the (unique) invariant control set of $S$ on $%
\mathbb{F}_{\Theta }$. It is closed, $S$-invariant and has nonempty
interior. Hence, it meets the dense set $\sigma _{\Theta }\left( H_{\mu
}\right) $.

\begin{lema}
If $G\left( \mu \right) \subset S$ then $C_{\Theta }\cap \mathrm{att}%
_{\Theta }\left( H_{\mu }\right) \neq \emptyset $.
\end{lema}

\begin{profe}
We have $C_{\Theta }\cap \sigma _{\Theta }\left( H_{\mu }\right) \neq
\emptyset $ and if $x\in C_{\Theta }\cap \sigma _{\Theta }\left( H_{\mu
}\right) $ then its $\omega $-limit $\omega \left( x\right) $ (w.r.t. $%
\widetilde{H}_{\mu }$) is contained in $C_{\Theta }$, because $\{\exp
tH_{\mu }:t\in \mathbb{R}\}\subset G\left( \mu \right) \subset S$. Since $%
\omega \left( x\right) \subset \mathrm{att}_{\Theta }\left( H_{\mu }\right) $%
, it follows that $\emptyset \neq \omega \left( x\right) \subset C_{\Theta
}\cap \mathrm{att}_{\Theta }\left( H_{\mu }\right) $.
\end{profe}

Now we look at the orbits $G\left( \mu \right) \cdot y$ through points $y\in 
\mathrm{att}_{\Theta }\left( H_{\mu }\right) =Z_{H_{\mu }}\cdot b_{\Theta }$%
. First for $y=b_{\Theta }$ we have the following general result.

\begin{lema}
\label{lemesse2}Let $\mathbb{F}_{\Theta }$ be a flag manifold and $\beta $ a
positive root. Then $G\left( \beta \right) \cdot b_{\Theta }$ is either a $2$%
-sphere or reduces to a point. If $\beta \notin \langle \Theta \rangle $
then $\dim G\left( \beta \right) \cdot b_{\Theta }=2$. In particular, $\dim
G\left( \mu \right) \cdot b_{\Theta }=2$ if $H_{\mu }\in \mathrm{cl}%
\mathfrak{h}_{\mathbb{R}}^{+}$.
\end{lema}

\begin{profe}
The point is that the orbit $G\left( \beta \right) \cdot b_{\Theta }$ equals 
$b_{\Theta }$ or identifies to the only flag manifold of $G\left( \beta
\right) $ which is the same as the flag manifold of $\mathrm{Sl}\left( 2,%
\mathbb{C}\right) $ (because $\mathfrak{g}\left( \beta \right) \approx 
\mathfrak{sl}\left( 2,\mathbb{C}\right) $), which in turn is $S^{2}$.

To see this denote by $\mathfrak{g}\left( \beta \right) _{b_{\Theta }}$ the
isotropy subalgebra at $b_{\Theta }$ for the action of $G\left( \beta
\right) $ on $\mathbb{F}_{\Theta }$. It contains the subalgebra $\mathfrak{p}%
_{\beta }=\mathrm{span}\{H_{\beta }\}\oplus \mathfrak{g}_{\beta }$ which is
a parabolic subalgebra of $\mathfrak{g}\left( \beta \right) $. This implies
that the isotropy subgroup at $b_{\Theta }$ contains the identity component
of the parabolic subgroup $P_{\beta }=\mathrm{Norm}_{G\left( \beta \right) }%
\mathfrak{p}_{\beta }\subset G\left( \beta \right) $. But any parabolic
subgroup of the complex group $G\left( \beta \right) $ is connected, hence $%
P_{\beta }$ is contained in the isotropy subgroup at $b_{\Theta }$, for the
action of $G$. This shows that $G\left( \beta \right) \cdot b_{\Theta }$ is
either a $2$-sphere or reduces to a point.

Now, if $\beta \notin \langle \Theta \rangle $ then $\mathfrak{g}_{-\beta }$
has zero intersection with the isotropy subalgebra $\mathfrak{p}_{\Theta }$
at $b_{\Theta }$ (which is the sum of the Cartan subalgra with root spaces).
This implies that $\mathfrak{g}\left( \beta \right) _{b_{\Theta }}=\mathfrak{%
p}_{\beta }$, and since an istropy subgroup normalizes the isotropy
subalgebra, it follows that $P_{\beta }$ is exactly the isotropy subgroup at 
$b_{\Theta }$ for the action of $G\left( \beta \right) $. Hence $G\left(
\beta \right) \cdot b_{\Theta }\approx G\left( \beta \right) /P_{\beta
}\approx S^{2}$.
\end{profe}

As to the $G\left( \mu \right) $-orbit through $y=g\cdot b_{\Theta }$, $g\in
Z_{H_{\mu }}$, we write 
\begin{equation*}
G\left( \mu \right) \cdot y=g\left( g^{-1}G\left( \mu \right) g\cdot
b_{\Theta }\right) 
\end{equation*}
so that $G\left( \mu \right) \cdot y$ is diffeomorphic to $g^{-1}G\left( \mu
\right) g\cdot b_{\Theta }$. The Lie group $g^{-1}G\left( \mu \right) g$ has
Lie algebra $\mathfrak{g}\left( \mu \right) ^{g}=\mathrm{Ad}\left( g\right)
\left( \mathfrak{g}\left( \mu \right) \right) $ also isomorphic to $%
\mathfrak{sl}\left( 2,\mathbb{C}\right) $. Since $\mathrm{Ad}\left( g\right)
H_{\mu }=H_{\mu }$ there is the root space decomposition 
\begin{equation*}
\mathfrak{g}\left( \mu \right) ^{g}=\langle H_{\mu }\rangle \oplus \mathrm{Ad%
}\left( g\right) \left( \mathfrak{g}_{\mu }\right) \oplus \mathrm{Ad}\left(
g\right) \left( \mathfrak{g}_{-\mu }\right) . 
\end{equation*}

\begin{lema}
Keep the assumption that $H_{\mu }\in \mathrm{cl}\mathfrak{h}_{\mathbb{R}%
}^{+}$ and take $g\in Z_{H_{\mu }}$. Then $\mathrm{Ad}\left( g\right) \left( 
\mathfrak{g}_{\mu }\right) \subset \mathfrak{n}^{+}=\sum_{\alpha >0}%
\mathfrak{g}_{\alpha }$.
\end{lema}

\begin{profe}
Since $\mathrm{Ad}\left( g\right) H_{\mu }=H_{\mu }$, it follows that $%
\mathrm{Ad}\left( g\right) $ commutes with $\mathrm{ad}\left( H_{\mu
}\right) $ and hence maps the eigenspaces of $\mathrm{ad}\left( H_{\mu
}\right) $ onto themselves. Now $\mathfrak{g}_{\mu }$ is contained in the $%
\mu \left( H_{\mu }\right) $-eigenspace of $\mathrm{ad}\left( H_{\mu
}\right) $. Hence, $\mathrm{Ad}\left( g\right) \mathfrak{g}_{\mu }$ is
contained in the same eigenspace. Now $\mu \left( H_{\mu }\right) =\langle
\mu ,\mu \rangle >0$ and the assumption that $H_{\mu }\in \mathfrak{h}_{%
\mathbb{R}}^{+}$ implies that the eigenspaces of $\mathrm{ad}\left( H_{\mu
}\right) $ associated to positive eigenvalues are contained in $\mathfrak{n}%
^{+}$. Hence $\mathrm{Ad}\left( g\right) \mathfrak{g}_{\mu }\subset 
\mathfrak{n}^{+}$ as claimed.
\end{profe}

%TCIMACRO{\TeXButton{vspace12pt}{\vspace{12pt}}}%
%BeginExpansion
\vspace{12pt}%
%EndExpansion

%TCIMACRO{\TeXButton{noindent}{\noindent}}%
%BeginExpansion
\noindent%
%EndExpansion
\textbf{Remark:} If $\mu $ is the highest root the above lemma has a more
precise statement, namely if $g\in Z_{H_{\mu }}$ then $g$ centralizes $%
\mathfrak{g}\left( \mu \right) $ and $G\left( \mu \right) $ (see Proposition %
\ref{propcentrisnorm} below). Hence $g^{-1}G\left( \mu \right) g=G\left( \mu
\right) $ and $G\left( \mu \right) \cdot y=g\left( G\left( \mu \right) \cdot
b_{\Theta }\right) $, what simplifies the proofs to follow.

\begin{lema}
Keep the above notation and assumptions. Then $G\left( \mu \right) \cdot
y=g\left( g^{-1}G\left( \mu \right) g\cdot b_{\Theta }\right) $ is either a $%
2$-sphere or reduces to a point.
\end{lema}

\begin{profe}
The subalgebra $\mathfrak{p}_{\mu }=\langle H_{\mu }\rangle \oplus \mathrm{Ad%
}\left( g\right) \mathfrak{g}_{\mu }$ is a parabolic subalgebra of $\mathrm{%
Ad}\left( g\right) \left( \mathfrak{g}\left( \mu \right) \right) $. Hence as
in the proof of Lemma \ref{lemesse2} it is enough to check that $\mathfrak{p}%
_{\mu }$ is contained in the isotropy subalgebra $\mathfrak{g}\left( \beta
\right) _{b_{\Theta }}^{g}$ (for the action of  $g^{-1}G\left( \mu \right) g$%
) at $b_{\Theta }$. Clearly $H_{\mu }\in \mathfrak{g}\left( \beta \right)
_{b_{\Theta }}^{g}$. On the other hand by the previous lemma $\mathrm{Ad}%
\left( g\right) \mathfrak{g}_{\mu }\subset \mathfrak{n}^{+}$ which in turn
is contained in the isotropy subalgebra at $b_{\Theta }$ (for the action of $%
G$). Hence $\mathfrak{p}_{\mu }\subset \mathfrak{g}\left( \mu \right)
_{b_{\Theta }}^{g}$ and the lemma follows.
\end{profe}

This proof shows that the orbit $g^{-1}G\left( \mu \right) g\cdot b_{\Theta
} $ induces a map from the flag manifold $G\left( \mu \right) /P_{\mu
}\approx S^{2}$ into $\mathbb{F}_{\Theta }$. We denote this map by $\sigma
_{g,\mu }:S^{2}\rightarrow \mathbb{F}_{\Theta }$.

The next, and final step, is to check that the $2$-spheres appearing in the
last lemma are not homotopic to a point (at least a great amount of them).

The idea is to exhibit a differential $2$-form $\Omega $ on $\mathbb{F}%
_{\Theta }$ with $d\Omega =0$ such that the pull-back $\nu =\sigma _{g,\mu
}^{\ast }\Omega $ is a (non zero) volume form on $S^{2}$. This would prove
that the map $\sigma _{g,\mu }^{\ast }:H^{2}\left( \mathbb{F}_{\Theta },%
\mathbb{R}\right) \rightarrow H^{2}\left( S^{2},\mathbb{R}\right) $ induced
on cohomology by $\sigma _{g,\mu }$ is not trivial, implying that $\sigma
_{g,\mu }$ is not homotopic to a constant map.

In fact, a volume form $\nu $ on the orientable manifold $S^{2}$ is a
generator of its $1$-dimensional de Rham cohomology $H^{2}\left( S^{2},%
\mathbb{R}\right) $, that is, $d\nu =0$ and $\nu $ is not $d\eta $ for a $1$%
-form $\eta $. If $\nu =\sigma _{g,\mu }^{\ast }\Omega $ then $\Omega $ is
not exact, for otherwise $\Omega =d\omega $ imply $\nu =\sigma _{\mu }^{\ast
}d\omega =d\sigma _{\mu }^{\ast }\omega $ and $\nu $ would be exact as well.
Hence $\Omega $ represents a non zero element in the De Rham cohomology $%
H^{2}\left( \mathbb{F}_{\Theta },\mathbb{R}\right) $ and the image of its
cohomology class under $\sigma _{g,\mu }^{\ast }$ is the cohomology class $%
[\nu ]\neq 0$.

A $2$-form $\Omega $ that does the job is a $U$-invariant symplectic form
associated to an invariant Hermitian metric together with a complex
structure on $\mathbb{F}_{\Theta }$ (K\"{a}hler form). The construction of
these geometric objects goes back to Borel \cite{bor}. To define it we
follow \cite{smher}. First we need a special basis of the tangent space $%
T_{b_{\Theta }}\mathbb{F}_{\Theta }$ at the origin. To get it start with a
Weyl basis of $\mathfrak{g}$ which is given by the choice of a generator $%
X_{\alpha }$ of the root space $\mathfrak{g}_{\alpha }$, for each root $%
\alpha $ and satisfying the conditions $\langle X_{\alpha },X_{-\alpha
}\rangle =1$ and $\left[ X_{\alpha },X_{\beta }\right] =m_{\alpha ,\beta
}X_{\alpha +\beta }$ with $m_{\alpha ,\beta }\in \mathbb{R}$ (see \cite%
{smher}, for details). Then define $A_{\alpha }=X_{\alpha }-X_{-\alpha }$, $%
S_{\alpha }=i\left( X_{\alpha }+X_{-\alpha }\right) $ and $\mathfrak{u}%
_{\alpha }=\mathrm{span}\{A_{\alpha },S_{\alpha }\}$ with $\alpha $ a
positive roots. Both $A_{\alpha }$ and $S_{\alpha }$ belong to the compact
real form $\mathfrak{u}$ of $\mathfrak{g}$ (Lie algebra of $U$). By the
action of $U$ on $\mathbb{F}_{\Theta }$ we get the induced vector fields $%
\widetilde{A}_{\alpha }$ and $\widetilde{S}_{\alpha }$. Then we have

\begin{itemize}
\item the set $\widetilde{\mathfrak{u}}_{\alpha }=\{\widetilde{A}_{\alpha
}\left( b_{\Theta }\right) ,\widetilde{S}_{\alpha }\left( b_{\Theta }\right)
:\alpha \notin \langle \Theta \rangle \}$ is a basis of $T_{b_{\Theta }}%
\mathbb{F}_{\Theta }$. \-
\end{itemize}

Now, $\Omega $ is defined by specifying its value $\Omega _{0}$ at the
origin and then extending to the whole $\mathbb{F}_{\Theta }$ by the action
of $U$. The extension is possible if $\Omega _{0}$ is invariant by the
isotropy representation of $U_{\Theta }$ on $T_{b_{\Theta }}\mathbb{F}%
_{\Theta }$.

To get $\Omega _{0}$ we choose first real numbers $\lambda _{\alpha }>0$, $%
\alpha \in \Pi ^{+}$, satisfying

\begin{enumerate}
\item $\lambda _{\alpha +\beta }=\lambda _{\alpha }+\lambda _{\beta }$ if $%
\alpha $, $\beta $ and $\alpha +\beta $ are positive roots.

\item $\lambda _{\alpha +\gamma }=\lambda _{\alpha }$ if $\alpha $, $\gamma $
and $\alpha +\gamma $ are roots with $\gamma \in \langle \Theta \rangle $.
\end{enumerate}

%TCIMACRO{\TeXButton{vspace12pt}{\vspace{12pt}}}%
%BeginExpansion
\vspace{12pt}%
%EndExpansion

%TCIMACRO{\TeXButton{noindent}{\noindent}}%
%BeginExpansion
\noindent%
%EndExpansion
\textbf{Remark:} Although it will not be used below we note that the numbers 
$\lambda _{\alpha }>0$, $\alpha \in \Pi ^{+}$, define an inner product on $%
T_{b_{\Theta }}\mathbb{F}_{\Theta }$, which by the second condition is $%
U_{\Theta }$-invariant, and hence extends to a Riemmannian metric $g$ on $%
\mathbb{F}_{\Theta }$. The first condition ensures that a Hermitian metric
built from $g$ and a complex structure $J$ on $\mathbb{F}_{\Theta }$ has a K%
\"{a}hler form which is symplectic (see \cite{smher}, Section 2.4).

Now we define $\Omega _{0}$ by declaring that $\Omega _{0}\left( X,Y\right)
=0$ if $X=\widetilde{\mathfrak{u}}_{\alpha }$ and $Y=\widetilde{\mathfrak{u}}
$ with $\alpha \neq \beta $ and 
\begin{equation*}
\Omega _{0}\left( \widetilde{A}_{\alpha }\left( b_{\Theta }\right) ,%
\widetilde{S}_{\alpha }\left( b_{\Theta }\right) \right) =\lambda _{\alpha
}\qquad \alpha \in \Pi ^{+}. 
\end{equation*}
The second condition above ensures that $\Omega _{0}$ is a $2$-form on $%
T_{b_{\Theta }}\mathbb{F}_{\Theta }$ invariant by $U_{\Theta }$, and hence
defines a $2$-form $\Omega $ on $\mathbb{F}_{\Theta }$, by translation. On
the other hand from the first condition we have $d\Omega =0$ (see \cite%
{smher}, Proposition 2.1).

Now we are prepared to prove that the $2$-spheres are not homotopic to a
point.

\begin{lema}
\label{lemdenso}Let $\mu $ be a positive root such that $H_{\mu }\in 
\mathfrak{h}_{\mathbb{R}}^{+}$, and denote by $Z_{H_{\mu }}$ the centralizer
of $H_{\mu }$ in $G$, and put $U_{H_{\mu }}=Z_{H_{\mu }}\cap U$. Take a flag
manifold $\mathbb{F}_{\Theta }$. Then there exists a subset $V\subset 
\mathrm{att}_{\Theta }\left( H_{\mu }\right) =U_{H_{\mu }}\cdot b_{\Theta }$
open and dense in $\mathrm{att}_{\Theta }\left( H_{\mu }\right) $ such that
for every $x\in V$, the orbit $G\left( \mu \right) \cdot x$ is a $2$-sphere
not homotopic to a point.
\end{lema}

\begin{profe}
For any $x\in \mathrm{att}_{\Theta }\left( H_{\mu }\right) $ we write $%
x=g\cdot b_{\Theta }$ with $g\in U_{H_{\mu }}$. Since $H_{\mu }\in \mathfrak{%
h}_{\mathbb{R}}^{+}$, we have by Proposition \ref{propsuport} that $\mathrm{%
supp}\mu =\Sigma $, so that by Lemma \ref{lemesse2}, $G\left( \mu \right)
\cdot b_{\Theta }$ is a $2$-sphere. Its tangent space $T_{b_{\Theta }}\left(
G\left( \mu \right) \cdot b_{\Theta }\right) $ has the basis $\{\widetilde{A}%
_{\mu }\left( b_{\Theta }\right) ,\widetilde{S}_{\mu }\left( b_{\Theta
}\right) \}$. Analogously, the tangent space at $b_{\Theta }$ of $%
g^{-1}G\left( \mu \right) g\cdot b_{\Theta }$, which we denote simply by $%
T^{g}$, is spanned by $\{\widetilde{A}_{\mu }^{g}\left( b_{\Theta }\right) ,%
\widetilde{S}_{\mu }^{g}\left( b_{\Theta }\right) \}$ where $A_{\mu }^{g}=%
\mathrm{Ad}\left( g\right) A_{\mu }$ and $S_{\mu }^{g}=\mathrm{Ad}\left(
g\right) S_{\mu }$. Either both vectors $\widetilde{A}_{\mu }^{g}\left(
b_{\Theta }\right) $ and $\widetilde{S}_{\mu }^{g}\left( b_{\Theta }\right) $
are zero or they form a basis of $T^{g}$.

Now we pull-back the symplectic form $\Omega $ to $T^{g}$ and define the
function $\phi :U_{H_{\mu }}\cdot b_{\Theta }\rightarrow \mathbb{R}$ by 
\begin{equation*}
\phi \left( g\cdot b_{\Theta }\right) =\Omega \left( \widetilde{A}_{\mu
}^{g}\left( b_{\Theta }\right) ,\widetilde{S}_{\mu }^{g}\left( b_{\Theta
}\right) \right) \qquad g\in U_{H_{\mu }}, 
\end{equation*}
which is well defined because any $g\in U_{H_{\mu }}$ leaves $\Omega $
invariant and $g\cdot b_{\Theta }=g_{1}\cdot b_{\Theta }$ implies that $\{%
\widetilde{A}_{\mu }^{g}\left( b_{\Theta }\right) ,\widetilde{S}_{\mu
}^{g}\left( b_{\Theta }\right) \}$ and $\{\widetilde{A}_{\mu }^{g_{1}}\left(
b_{\Theta }\right) ,\widetilde{S}_{\mu }^{g_{1}}\left( b_{\Theta }\right) \}$
span the same subspace, namely $T^{g}=T^{g_{1}}$. Hence $\Omega \left( 
\widetilde{A}_{\mu }^{g}\left( b_{\Theta }\right) ,\widetilde{S}_{\beta
}^{g}\left( b_{\Theta }\right) \right) =\Omega \left( \widetilde{A}_{\mu
}^{g_{1}}\left( b_{\Theta }\right) ,\widetilde{S}_{\mu }^{g_{1}}\left(
b_{\Theta }\right) \right) $

The function $\phi $ is analytic as is the map $g\mapsto \mathrm{Ad}\left(
g\right) $. It is not idencally zero, since by Lemma \ref{lemesse2} we have 
\begin{equation*}
\phi \left( b_{\Theta }\right) =\Omega \left( \widetilde{A}_{\mu }\left(
b_{\Theta }\right) ,\widetilde{S}_{\mu }\left( b_{\Theta }\right) \right)
=\lambda _{\mu }\neq 0. 
\end{equation*}

Hence the subset $V=\{x\in U_{H_{\mu }}\cdot b_{\Theta }:\phi \left(
x\right) \neq 0\}$ is open and dense in $U_{H_{\mu }}\cdot b_{\Theta }$. For
any $x\in V$ the orbit $G\left( \mu \right) \cdot x$ is $2$-dimensional.
Also, if $x=g\cdot b_{\Theta }$, $g\in U_{H_{\mu }}$, then $\Omega \left( 
\widetilde{A}_{\mu }^{g}\left( b_{\Theta }\right) ,\widetilde{S}_{\mu
}^{g}\left( b_{\Theta }\right) \right) \neq 0$. Now $\Omega $ is invariant
by $\left( g^{-1}G\left( \mu \right) g\right) \cap U$ and this group acts
transitively on $g^{-1}G\left( \mu \right) g\cdot b_{\Theta }$. Hence the
pull-back of $\Omega $ to $g^{-1}G\left( \mu \right) g\cdot b_{\Theta }$ is
a volume form, which shows that $S^{2}\approx g^{-1}G\left( \mu \right)
g\cdot b_{\Theta }$ is not homotopic to a point if $x=g\cdot b_{\Theta }\in
V $.

In conclusion, we have $G\left( \mu \right) \cdot x=g\left( g^{-1}G\left(
\mu \right) g\cdot b_{\Theta }\right) $, so that $G\left( \mu \right) \cdot
x $ is not homotopic to a point as well.
\end{profe}

%TCIMACRO{\TeXButton{vspace12pt}{\vspace{12pt}}}%
%BeginExpansion
\vspace{12pt}%
%EndExpansion

%TCIMACRO{\TeXButton{noindent}{\noindent}}%
%BeginExpansion
\noindent%
%EndExpansion
\textbf{End of proof of Theorem }\ref{teocomplexo}: If $x=ng\cdot b_{\Theta
}\in \sigma _{\Theta }\left( H_{\mu }\right) $ with $n\in N^{-}$ and $g\in
Z_{H_{\mu }}$ then $\left( \exp tH_{\mu }\right) ng\cdot b_{\Theta }=\left(
\exp tH_{\mu }\right) n\left( \exp \left( -tH_{\mu }\right) \right) g\cdot
b_{\Theta }$ because $\left( \exp \left( -tH_{\mu }\right) \right) g=g\left(
\exp \left( -tH_{\mu }\right) \right) $ and $\left( \exp \left( -tH_{\mu
}\right) \right) \cdot b_{\Theta }=b_{\Theta }$. However $\lim_{t\rightarrow
\infty }\left( \exp tH_{\mu }\right) n\left( \exp \left( -tH_{\mu }\right)
\right) =1$, so that 
\begin{equation*}
\lim_{t\rightarrow +\infty }\left( \exp tH_{\mu }\right) ng\cdot b_{\Theta
}=g\cdot b_{\Theta }.
\end{equation*}%
Now, let $V$ be the open and dense subset of $Z_{H_{\beta }}\cdot b_{\Theta }
$ ensured by the last lemma. Then $N_{\mu }^{-}\cdot V$ is open and dense in 
$\mathbb{F}_{\Theta }$. Since $C_{\Theta }$ has non empty interior we can
find $x\in C_{\Theta }$ such that 
\begin{equation*}
\lim_{t\rightarrow +\infty }\left( \exp tH\right) x=y\in C_{\Theta }\cap V.
\end{equation*}%
Then $G\left( \mu \right) \cdot y$ is a $2$-sphere not homotopic to a point
contained in $C_{\Theta }$. This shows that $C_{\Theta }$ cannot be
contained in a contractible subset of $\mathbb{F}_{\Theta }$. Since $\Theta $
was arbitrary $S=G$. In view of Lemma \ref{lemactweyl}, this proves Theorem %
\ref{teocomplexo}.  
%TCIMACRO{\TeXButton{qed}{\qed}}%
%BeginExpansion
\qed%
%EndExpansion

To conclude this section we prove the following statement ensuring that for
the highest root we have $gG\left( \mu \right) g^{-1}=G\left( \mu \right) $, 
$g\in Z_{H_{\mu }}$, so that the set $V$ of Lemma \ref{lemdenso} is the
totality of $\mathrm{att}_{\Theta }\left( H_{\mu }\right) $.

\begin{proposicao}
\label{propcentrisnorm}Let $\mu $ be the highest root, and suppose that $%
g\in G$ centralizes $H_{\mu }$ that is $\mathrm{Ad}\left( g\right) H_{\mu
}=H_{\mu }$. Then $g$ normalizes $G\left( \mu \right) $ (actually $g$
comutes with every $h\in G\left( \mu \right) $).
\end{proposicao}

\begin{profe}
Write $Z_{H_{\mu }}$ for the centralizer of $H_{\mu }$ in $G$. It is a Lie
group with Lie algebra $\mathfrak{z}_{H_{\mu }}=\ker \mathrm{ad}\left(
H_{\mu }\right) $, the centralizer of $H_{\mu }$ in $\mathfrak{g}$. By the
root space decomposition we have 
\begin{equation*}
\mathfrak{z}_{H_{\mu }}=\mathfrak{h}+\sum_{\beta \left( H_{\mu }\right) =0}%
\mathfrak{g}_{\beta }. 
\end{equation*}
Write $\Theta _{H_{\mu }}=\{\alpha \in \Sigma :\alpha \left( H_{\mu }\right)
=0\}$. Then a root $\beta $ anihilates $H_{\mu }$ if and only if $\mathrm{%
supp}\beta \subset \Theta _{H_{\mu }}$. This follows from the fact that $%
H_{\mu }\in \mathrm{cl}\mathfrak{h}_{\mathbb{R}}^{+}$, so that if $\alpha
\notin \Theta _{H_{\mu }}$ then $\alpha \left( H_{\mu }\right) >0$.

Take a root $\beta $ with $\beta \left( H_{\mu }\right) =\langle \beta ,\mu
\rangle =0$. Then $\mu \pm \beta $ are not roots. In fact if $\beta >0$ then 
$\mu +\beta $ is not a root since $\beta >0$. Hence by the Killing formula,
the orthogonality $\langle \beta ,\mu \rangle =0$ implies that $\mu -\beta $
is neither a root. We have also that $-\mu \pm \beta $ are not roots.
Therefore $[\mathfrak{g}_{\pm \mu },\mathfrak{g}_{\beta }]=0$ if $\beta
\left( H_{\mu }\right) =0$, and since $\mathfrak{h}$ is abelian we conclude
that $[X,\mathfrak{g}\left( \mu \right) ]=0$ if $X\in \mathfrak{z}_{H_{\mu
}} $.

Now, since we are working with the complex group $G$ it is true that $%
Z_{H_{\mu }}$ is connected. Hence the commutativity between $\mathfrak{z}%
_{H_{\mu }}$ and $\mathfrak{g}\left( \mu \right) $ implies that $\mathrm{Ad}%
\left( g\right) Y=Y$ for any $Y\in \mathfrak{z}_{H_{\mu }}$. This in turn
implies the elements of $Z_{H_{\mu }}$ commute with the elements of $G\left(
\mu \right) $.
\end{profe}

\begin{corolario}
\label{cornopoint}Let $\mu $ be the highest root and denote by $Z_{H_{\mu }}$
the centralizer of $H_{\mu }$ in $G$. Take a flag manifold $\mathbb{F}%
_{\Theta }$. Then for any $g\in Z_{H_{\mu }}$ the orbit $G\left( \mu \right)
\cdot gb_{\Theta }$ is a $2$-sphere in $\mathbb{F}_{\Theta }$ not homotopic
to a point.
\end{corolario}

\section{Controllability theorem}

As mentioned in the introduction our source of inspiration for Theorem \ref%
{teocomplexo} are the results on controllability of control systems of \cite%
{JK}, \cite{JK1}, \cite{GKS}, \cite{EGK}. The starting point in the proof of
these results is the proof that $G\left( \mu \right) $ is contained in the
semigroup of control. Their assumptions are designed to ensure this
inclusion. With Theorem \ref{teocomplexo} we can improve (for complex Lie
groups) the final theorem of \cite{EGK}, without insisting to work with the
highest root.

Let 
\begin{equation}
\dot{g}=\left( A+u\left( t\right) B\right) g\qquad u\left( t\right) \in 
\mathbb{R}  \label{forcontrol}
\end{equation}
be a rigth invariant control system with unrestricted controls where $A$,$%
B\in \mathfrak{g}$ with $\mathfrak{g}$ a complex Lie algebra of the complex
Lie group $G$. We let $S$ be the semigroup of the system (generated by $\exp
t\left( A+uB\right) $, $t\geq 0$, $u\in \mathbb{R}$) and denote by 
\begin{equation*}
\Gamma =\{X\in \mathfrak{g}:\forall t\geq 0,\,\exp tX\in \mathrm{cl}S\} 
\end{equation*}
its Lie wedge. $\Gamma $ is a closed convex cone invariant by $\exp t\mathrm{%
ad}\left( X\right) $, $t\in \mathbb{R}$, if $\pm X\in \Gamma $ (see
Hilgert-Hofmann-Lawson \cite{hhl}).

Since (\ref{forcontrol}) is with unrestricted controls the following easy
argument shows that $\pm B\in \Gamma $: $A+uB\in \Gamma $, $u\in \mathbb{R}$%
, hence if $u>0$, $\left( 1/u\right) A+B\in \Gamma $, so that $%
B=\lim_{u\rightarrow +\infty }\left( \left( 1/u\right) A+B\right) \in \Gamma 
$. Similarly $-B\in \Gamma $, by making $u\rightarrow -\infty $.

Now we shall take $B$ in the Cartan subalgebra $\mathfrak{h}$ and write 
\begin{equation*}
A=A_{0}+\sum_{\alpha \in \Pi }A_{\alpha } 
\end{equation*}
for the root space decomposition of $A$, $A_{0}\in \mathfrak{h}$ and $%
A_{\alpha }\in \mathfrak{g}_{\alpha }$.

The Cartan subalgebra $\mathfrak{h}$ decomposes as $\mathfrak{h}=\mathfrak{h}%
_{\mathbb{R}}+i\mathfrak{h}_{\mathbb{R}}$ where $\mathfrak{h}_{\mathbb{R}}$
is the real subspace where the roots assume real values. If $\beta $ is a
root we have $\beta \left( H\right) \in \mathbb{R}$ if $H\in \mathfrak{h}_{%
\mathbb{R}}$ and $\beta \left( H\right) $ is immaginary if $H\in i\mathfrak{h%
}_{\mathbb{R}}$.

In particular we write $B=B_{\mathrm{Re}}+B_{\mathrm{Im}}\in \mathfrak{h}_{%
\mathbb{R}}+i\mathfrak{h}_{\mathbb{R}}$, and state the controllability
result separetaly into two cases: 1) $\mathrm{ad}\left( B\right) $ has
purely imaginary eigenvalues, that is, $B_{\mathrm{Re}}=0$; 2) $B_{\mathrm{Re}%
}\neq 0$. The proofs follow almost immediately from our Theorem \ref%
{teocomplexo} and Lemma 2.3 of \cite{EGK}, whose arguments we reproduce, for
the sake of completeness.

\begin{teorema}
In the control system (\ref{forcontrol}) take $B\in \mathfrak{h}$ and
suppose that $\mathrm{ad}\left( B\right) $ has purely imaginary eigenvalues (%
$B_{\mathrm{Re}}=0$). Then the system is controllable in $G$ if

\begin{enumerate}
\item $A$ and $B$ generate $\mathfrak{g}$ (Lie algebra rank condition), and

\item there exists a root $\alpha $ such that $\alpha \left( B\right) \neq 0 
$, $A_{\alpha }\neq 0\neq A_{-\alpha }$ and $\alpha \left( B\right) \neq
\beta \left( B\right) $ for any root $\beta \neq \alpha $ with $A_{\beta
}\neq 0$.
\end{enumerate}
\end{teorema}

\begin{profe}
We have $\exp t\mathrm{ad}\left( B\right) A\in \Gamma $ for all $t\in 
\mathbb{R}$ and 
\begin{equation*}
\exp t\mathrm{ad}\left( B\right) A=\sum_{\beta \in \Pi }e^{t\beta \left(
B\right) }A_{\beta }. 
\end{equation*}
Hence $A_{\eta }\left( t\right) =\left( 1+\eta \cos \alpha \left( B\right)
\right) \exp t\mathrm{ad}\left( B\right) A\in \Gamma $ if $|\eta |<1$. Since 
$\beta \left( B\right) \neq \alpha \left( B\right) $ is purely immaginary we
have 
\begin{equation*}
\lim_{T\rightarrow +\infty }\left( 1/T\right) \int_{0}^{T}e^{t\beta \left(
B\right) }\left( 1+\eta \cos \alpha \left( B\right) \right) dt=0, 
\end{equation*}
yielding the limit $\lim_{T\rightarrow +\infty }\left( 1/T\right)
\int_{0}^{T}A_{\eta }\left( t\right) =\left( \eta /2\right) A_{\alpha }$
(see \cite{JK1} and \cite{EGK}, Lemma 2.3). Therefore $A_{\alpha }\neq 0$
belongs to $\Gamma $, hence $\exp t\mathrm{ad}\left( B\right) A_{\alpha }\in
\Gamma $, $t\in \mathbb{R}$. Now, $\exp t\mathrm{ad}\left( B\right)
A_{\alpha }=e^{t\alpha \left( B\right) }A_{\alpha }$, and since $\alpha
\left( B\right) \neq 0$ we see that the complex subspace spanned by $%
A_{\alpha }$ is contained in $\Gamma $, that is, $\mathfrak{g}_{\alpha
}\subset \Gamma $. The same way it follows that $\mathfrak{g}_{-\alpha
}\subset \Gamma $, therefore $\mathfrak{g}\left( \alpha \right) \subset
\Gamma $ implying that $G\left( \alpha \right) \subset \Gamma $ and $S=G$,
by Theorem \ref{teocomplexo}.
\end{profe}

%TCIMACRO{\TeXButton{vspace12pt}{\vspace{12pt}}}%
%BeginExpansion
\vspace{12pt}%
%EndExpansion

%TCIMACRO{\TeXButton{noindent}{\noindent}}%
%BeginExpansion
\noindent%
%EndExpansion
\textbf{Remark:} In \cite{EGK} and \cite{JK1} the above result is proved
with the assumption that $B$ is strong regular, which means that $\alpha
\left( B\right) \neq 0$ for any root $\alpha $ and $\alpha \left( B\right)
\neq \beta \left( B\right) $ for roots $\alpha \neq \beta $. With strong
regularity it is possible to prove that $\mathfrak{g}\left( \alpha \right)
\subset \Gamma $ for several roots $\alpha $  and conclude that $\Gamma =%
\mathfrak{g}$. By applying Theorem \ref{teocomplexo} it is enough to have $%
\mathfrak{g}\left( \alpha \right) \subset \Gamma $ for just one root $\alpha 
$.

\begin{teorema}
\label{teocontrolreal}In the control system (\ref{forcontrol}) take $B\in 
\mathfrak{h}$ with $B_{\mathrm{Re}}\neq 0$. Then the system is controllable in 
$G$ if $A$ and $B$ generate $\mathfrak{g}$ (Lie algebra rank condition) and
there exists a root $\alpha $ such that

\begin{enumerate}
\item $\mathrm{Im}\alpha \left( B\right) \neq 0$.

\item If $\beta \neq \alpha $ is a positive root such that $\mathrm{Re}\beta
\left( B\right) \leq \mathrm{Re}\alpha \left( B\right) $ then $\mathrm{Re}\beta
\left( B\right) <\mathrm{Re}\alpha \left( B\right) $.

\item $A_{\pm \alpha }\neq 0$ and $A_{\gamma }=0$ in case $\mathrm{Re}\gamma
\left( B\right) >\mathrm{Re}\alpha \left( B\right) $ or $\mathrm{Re}\gamma
\left( B\right) <-\mathrm{Re}\alpha \left( B\right) $.
\end{enumerate}
\end{teorema}

\begin{profe}
For all $t\in \mathbb{R}$, $e^{\pm t\mathrm{Re}\left( \alpha \left( B\right)
\right) }\exp t\mathrm{ad}\left( B\right) A\in \Gamma $. By the third
condition 
\begin{equation*}
\exp t\mathrm{ad}\left( B\right) A=\sum e^{t\gamma \left( B\right) }A_{\beta
}
\end{equation*}%
with the sum extended to $\gamma $ with -$\alpha \left( B\right) \leq \gamma
\left( B\right) \leq \alpha \left( B\right) $. But by the second condition 
\begin{equation*}
\lim_{t\rightarrow +\infty }e^{-t\mathrm{Re}\alpha \left( B\right) }\exp t%
\mathrm{ad}\left( B\right) A=A_{\alpha }\quad \mathrm{and}\quad
\lim_{t\rightarrow -\infty }e^{t\mathrm{Re}\alpha \left( B\right) }\exp t%
\mathrm{ad}\left( B\right) A=A_{-\alpha },
\end{equation*}%
hence $A_{\pm \alpha }\in \Gamma $. Now, $\exp t\mathrm{ad}\left( B\right)
A_{\pm \alpha }=e^{\pm t\alpha \left( B\right) }A_{\pm \alpha }$, and since $%
\alpha \left( B\right) \neq 0$ we conclude that $\mathfrak{g}_{\pm \alpha
}\subset \Gamma $. Hence $\mathfrak{g}\left( \alpha \right) \subset \Gamma $%
, $G\left( \alpha \right) \subset S$ and the result follows by Theorem \ref%
{teocomplexo}.
\end{profe}

%TCIMACRO{\TeXButton{vspace12pt}{\vspace{12pt}}}%
%BeginExpansion
\vspace{12pt}%
%EndExpansion

%TCIMACRO{\TeXButton{noindent}{\noindent}}%
%BeginExpansion
\noindent%
%EndExpansion
\textbf{Remark:} In \cite{EGK} and \cite{JK1} the above theorem is proved by
taking the  highest root instead of an arbitrary root $\alpha $. In fact, if 
$B_{\mathrm{Re}}\in \mathfrak{h}_{\mathbb{R}}^{+}$ (which can be assumed
without loss of generality) then the second condition and part of the third
condition are automatically true when $\alpha $ is the corresponding highest
root. In this case the assumption in \cite{EGK} and \cite{JK1} is that $%
A_{\pm \alpha }\neq 0$. As to the first condition it follows if $B$ is
strong regular in the sense of \cite{EGK} and \cite{JK1}. This means that
the dimension of $\ker \mathrm{ad}\left( B\right) $ is the rank of $%
\mathfrak{g}$ and the eigenvalues of the complexification $\mathrm{ad}\left(
B\right) _{\mathbb{C}}$ of $\mathrm{ad}\left( B\right) $ are simple. In the
complex Lie algebra $\mathfrak{g}$ we must complexify its realification.
Then the eigenvalues of $\mathrm{ad}\left( B\right) _{\mathbb{C}}$ are those
of $\mathrm{ad}\left( B\right) $ together with their complex conjugates.
Hence the eigenvalues of $\mathrm{ad}\left( B\right) _{\mathbb{C}}$ are
simple if and only if no eigenvalue of $\mathrm{ad}\left( B\right) $ is
real. Therefore the strong regular condition implies that $\mathrm{Im}\beta
\left( B\right) \neq 0$ for any root $\beta $.


\begin{thebibliography}{99}
\bibitem{EGK} R. El Assoudi, J. P. Gauthier and I. Kupka,\textit{\ On
subsemigroups of semisimple Lie groups.} Annales de l'H. P., section $3$, 
\textbf{13} (1996), 117-133.

\bibitem{bor} A. Borel, \textit{K\"{a}hlerian coset spaces of semi-simple
Lie groups}. Proc. Nat. Acad. Sci. \textbf{40} (1954), 1147-1151.

\bibitem{smbrsimplt} C. J. Braga Barros and L.A.B San Martin, \textit{%
Controllability of Discrete-time Control Systems on the Symplectic Group}.
Systems \& Control Letters \textbf{42} (2001), 95-100.

\bibitem{dkv} J.J. Duistermat, J.A.C. Kolk and V.S. Varadarajan, \ \textit{\
Functions, flows and oscilatory integral on flag manifolds}. Compos.\ Math.\
49, 309-398 (1983).

\bibitem{fepase} T. Ferraiol, M. Patr\~{a}o and L. Seco, \textit{Jordan
decomposition and dynamics on flag manifolds}. Discrete and Continuous
Dynamical Systems, v. 26, p. 923-947, 2010.

\bibitem{GKS} J.P. Gauthier, I. Kupka and G. Sallet, \textit{Controllability
of simple Lie groups right invariant systems on real.} Systems Control
Letters \textbf{5} (1984) 187-190.

\bibitem{hhl} {Hilgert, J., K. H. Hofmann, and J. D. Lawson :} \textrm{``Lie
Groups, Convex Cones and Semigroups,''}\textrm{Oxford University Press, 1989}

\bibitem{JK} V. Jurdjevic, and I. Kupka, \textit{Control systems subordinate
to a group action: accessibility.} J. of Diff. Eq., \textbf{3}9 (1981),
186-211.

\bibitem{JK1} V.Jurdjevic, and I. Kupka, \textit{Control systems on
semisimple Lie groups and their homogeneous spaces.} Ann.Inst. Fourier
(Grenoble), \textbf{31}(1981), 151-179.

\bibitem{K} A.W Knapp, Lie groups beyond and introduction. Second edition,
Birkhauser (2004).

\bibitem{smroconex} O.G. do Rocio and L.A.B San Martin, \textit{\ Connected
components of open semigroups in semi-simple Lie groups}. Semigroup Forum 
\textbf{69} (2004), 1--29.

\bibitem{sminv} L.A.B San Martin, \textit{Invariant Control Sets on Flag
Manifolds}. Math. of Control, Signal and Systems, vol. 6 (1993), 41-61.

\bibitem{smsurv} L.A.B San Martin, \textit{Control Sets and Semigroups in
Semi-Simple Lie Groups}. In Semigroups in Algebra, Analysis and Geometry. De
Gruyter Expositions in Mathematics, vol. \textbf{20} (Editors: Hofmann, D.
H. , Lawson, J. e Vinberg, E. B.) (1995), 275-291.

\bibitem{SM} L.A.B San Martin, \textit{\ On global controllability of
discrete-time control systems.} Math. Control Signals Systems, \textbf{8}
(1995), 279-297.

\bibitem{SMT} L.A.B San Martin and P.A Tonelli, \textit{\ Semigroup actions
on homogeneous spaces.} Semigroups Forum, \textbf{50} (1995), 59-88.

\bibitem{smord} L.A.B San Martin, \textit{Order and Domains of Attraction of
Control Sets in Flag Manifolds}. Journal of Lie Theory, \textbf{8} (1998),
335-350.

\bibitem{smmax} L.A.B San Martin, \textit{Maximal semigroups in semi-simple
Lie groups}. Trans. Amer. Math. Soc., \textbf{353} (2001), 5165-5184.

\bibitem{smsanhom} L.A.B San Martin and A. J. Santana, \textit{Homotopy type
of Lie semigroups in semi-simple Lie groups}. Monatshefte f\"{u}r
Mathematik, \textbf{136} (2002), 151-173.

\bibitem{smher} L.A.B. San Martin and C.J.C. Negreiros, \textit{Invariant
almost Hermitian structures on flag manifolds}. Adv. Math., \textbf{178}
(2003), 277-310.

\bibitem{W} G. Warner, \textit{Harmonic Analysis on Semi-simple Lie Groups.}
Springer-Verlag, Berlin, (1972).
\end{thebibliography}
\end{document}